# A complete solution of the $k$-uniform supertrees with the eight largest $\alpha$-spectral radii

Lou-Jun Yu, Wen-Huan Wang[*]

*Department of Mathematics, Shanghai University, Shanghai 200444, China*

(Dated: February 12, 2023)

Let $\mathcal{T}(n,k)$ be the set of the $k$-uniform supertrees with $n$ vertices and $m$ edges, where $k \geq 3$, $n \geq 5$ and $m = \frac{n-1}{k-1}$. A conjecture concerning the supertrees with the fourth through the eighth largest $\alpha$-spectral radii in $\mathcal{T}(n,k)$ was proposed by You et al. (2020), where $0 \leq \alpha < 1$, $k \geq 3$ and $m \geq 10$. This conjecture was partially solved for $1 - \frac{1}{m-2} \leq \alpha < 1$ and $m \geq 10$ by Wang et al. (2022). When $0 \leq \alpha < 1 - \frac{1}{m-2}$ and $m \geq 10$, whether this conjecture is correct or not remains a problem to be further solved. By using a new $\rho_\alpha$-normal labeling method proposed in this article for computing the $\alpha$-spectral radius of the $k$-uniform hypergraphs, we completely prove that this conjecture is right for $0 \leq \alpha < 1$ and $m \geq 13$.

## 1. INTRODUCTION

Let $\mathcal{G} = (V(\mathcal{G}), E(\mathcal{G}))$ be a hypergraph, where $V(\mathcal{G}) = \{v_1, \ldots, v_n\}$ and $E(\mathcal{G}) = \{e_1, \ldots, e_m\}$ are the sets of the vertices and the edges of $\mathcal{G}$, respectively. If each edge of $\mathcal{G}$ has $k$ vertices, then $\mathcal{G}$ is a $k$-uniform hypergraph, where $k \geq 2$. Cooper and Dutle [1] defined the adjacency tensor of a $k$-uniform hypergraph $\mathcal{G}$ on $n$ vertices as the $k$-ordered and $n$-dimensional tensor $\boldsymbol{\mathcal{A}}(\mathcal{G}) = (a_{i_1 i_2 \cdots i_k})$, where $a_{i_1 i_2 \cdots i_k} = \frac{1}{(k-1)!}$ if $\{v_{i_1}, v_{i_2}, \ldots, v_{i_k}\} \in E(\mathcal{G})$ and $a_{i_1 i_2 \cdots i_k} = 0$ otherwise. Let $d_{\mathcal{G}}(v_i)$ be the degree of $v_i$ with $1 \leq i \leq n$. The degree diagonal tensor of $\mathcal{G}$ is the $k$-ordered and $n$-dimensional tensor $\boldsymbol{\mathcal{D}}(\mathcal{G}) = (d_{i_1 i_2 \cdots i_k})$, where $d_{i_1 i_2 \cdots i_k} = d_{\mathcal{G}}(v_i)$ if $i_j = i$ for $1 \leq j \leq k$ and $1 \leq i \leq n$, and $d_{i_1 i_2 \cdots i_k} = 0$ otherwise.

Inspired by the work of Nikiforov [2], Lin et al. [3] proposed to study the convex linear combinations $\boldsymbol{\mathcal{A}}_\alpha(\mathcal{G}) = \alpha \boldsymbol{\mathcal{D}}(\mathcal{G}) + (1-\alpha)\boldsymbol{\mathcal{A}}(\mathcal{G})$, where $0 \leq \alpha < 1$. The $\alpha$-spectral radius of $\mathcal{G}$, denoted by $\rho_\alpha(\mathcal{G})$, is defined to be the largest modulus of all the eigenvalues of $\mathcal{A}_\alpha(\mathcal{G})$. Obviously, when $\alpha = 0$, $\mathcal{A}_\alpha(\mathcal{G})$ is $\mathcal{A}(\mathcal{G})$ and $\rho_0(\mathcal{G})$ is the spectral radius of $\mathcal{G}$. When $\alpha = \frac{1}{2}$, $2\mathcal{A}_\alpha(\mathcal{G})$ is the signless Laplacian tensor of $\mathcal{G}$ and $2\rho_{\frac{1}{2}}(\mathcal{G})$ is the signless Laplacian spectral radius of $\mathcal{G}$.

---
[*] Corresponding author. Email: whwang@shu.edu.cn

The hypergraphs with the largest $\alpha$-spectral radii have been obtained among several classes of hypergraphs. For example, the $k$-uniform non-caterpillar hypergraphs with a given diameter [4], the hypergraphs with a given number of pendent edges [5], the unicyclic hypergraphs [5], the $k$-uniform unicyclic hypergraphs with a fixed diameter [6], and the $k$-uniform unicyclic hypergraphs with a given number of pendent edges [6], etc. For the upper bounds of the $\alpha$-spectral radius of hypergraphs, one can refer to Refs. [3, 5, 7–9].

Let $S_{m+1}$ be a star of order $m+1$, where $m \geq 2$. Let $S_{a,b}$ be a tree of order $m+1$ obtained from an edge $e = u_1 u_2$ by attaching $a$ and $b$ pendent edges to $u_1$ and $u_2$, respectively, where $a, b \geq 1$ with $a + b = m - 1$. Obviously, $S_{m+1} = S_{0,m-1}$. Let $S_{s_1,s_2,s_3}$ be a tree of order $m+1$ obtained from a path $P_3 = u_1 u_2 u_3$ by attaching $s_1$, $s_2$ and $s_3$ pendent edges to $u_1$, $u_2$ and $u_3$, respectively, where $s_1$, $s_2$ and $s_3$ are integers satisfying $s_1, s_3 \geq 1$, $s_2 \geq 0$ and $s_1 + s_2 + s_3 = m - 2$. The $k$-th power of a graph $H$, denoted by $\mathcal{H}^k$, is obtained from $H$ by adding $k - 2$ new vertices into each edge of $H$, where $k \geq 3$. A hypertree is the $k$-th power of an ordinary tree. Let $\mathcal{S}_{m+1}^k$, $\mathcal{S}_{a,b}^k$, and $\mathcal{S}_{s_1,s_2,s_3}^k$ be the $k$-th powers of $S_{m+1}$, $S_{a,b}$, and $S_{s_1,s_2,s_3}$, respectively. For example, $\mathcal{S}_{a,b}^k$ and $\mathcal{S}_{s_1,s_2,s_3}^k$ are shown in Figs. 1(a) and 1(b).

A supertree is a connected and acyclic hypergraph. Let $\mathcal{T}(t_1, t_2, t_3)$ be a $k$-uniform supertree obtained from $e_1 = \{u_1, u_2, \cdots, u_k\}$ by attaching $t_1$, $t_2$, and $t_3$ pendent edges having $k$ vertices at the vertices $u_1$, $u_2$, and $u_3$ of $e_1$, respectively, where $t_1, t_2, t_3 \geq 1$ with $t_1 + t_2 + t_3 = m - 1$ and $m \geq 4$. For example, $\mathcal{T}(t_1, t_2, t_3)$ is shown in Fig. 1(c).

Let $\mathcal{T}(n, k)$ be the set of the $k$-uniform supertrees of order $n$, where $k \geq 3$. For any supertree $\mathcal{T} \in \mathcal{T}(n, k)$, the number of the edges of $\mathcal{T}$ is $m = \frac{n-1}{k-1}$. For simplicity, let $\mathcal{T}^{(8)}(n, k)$ be the set of the eight preceding $k$-uniform supertrees in terms of the decreasing order of their spectral radii among $\mathcal{T}(n, k)$ [10, 11], where

$$\mathcal{T}^{(8)}(n, k) = \{\mathcal{S}_{m+1}^k, \mathcal{S}_{1,m-2}^k, \mathcal{S}_{2,m-3}^k, \mathcal{T}(1, 1, m-3), \mathcal{S}_{1,m-4,1}^k, \mathcal{S}_{m-3,0,1}^k, \mathcal{S}_{3,m-4}^k, \mathcal{T}(1, 2, m-4)\}.$$

Among $\mathcal{T}(n, k)$, You et al. [12] obtained that $\mathcal{S}_{m+1}^k$, $\mathcal{S}_{1,m-2}^k$, and $\mathcal{S}_{2,m-3}^k$ are the supertrees with the first, the second, and the third largest $\alpha$-spectral radii, respectively, where $k \geq 3$ and $m = \frac{n-1}{k-1} \geq 5$. For the supertrees with the fourth through the eighth

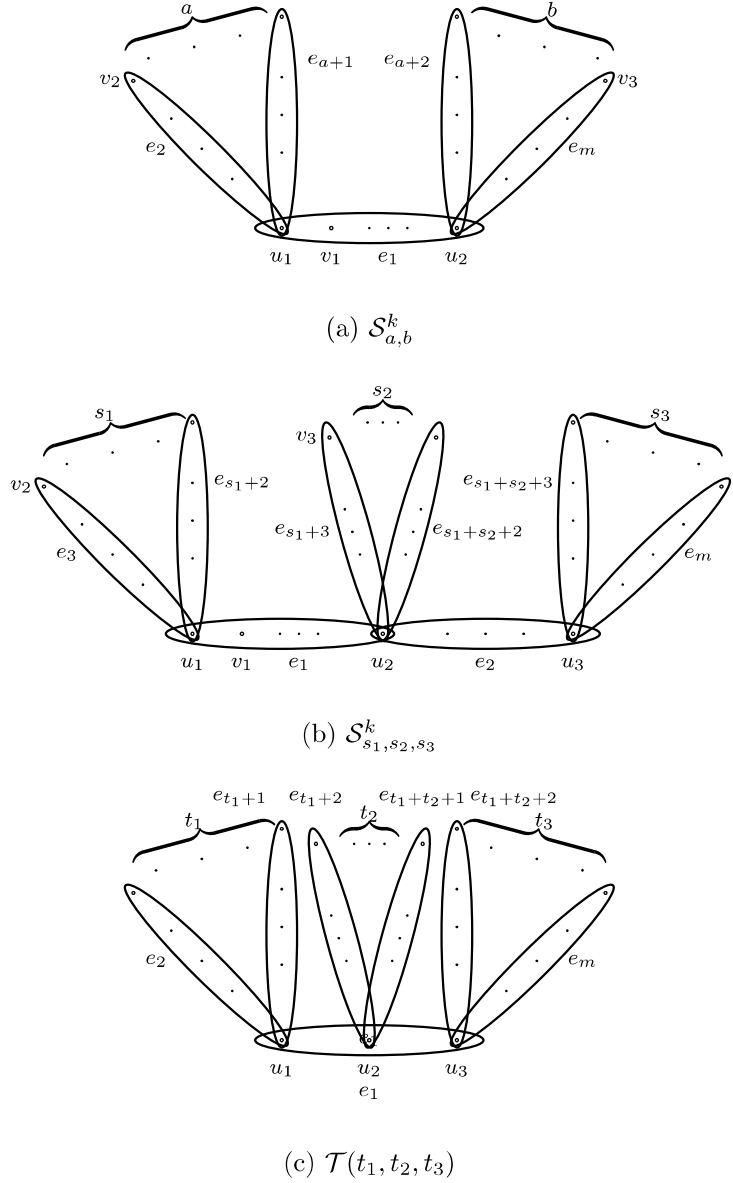

(a) $\mathcal{S}_{a,b}^k$

(b) $\mathcal{S}_{s_1,s_2,s_3}^k$

(c) $\mathcal{T}(t_1,t_2,t_3)$

Fig. 1: $\mathcal{S}_{a,b}^k$, $\mathcal{S}_{s_1,s_2,s_3}^k$ and $\mathcal{T}(t_1,t_2,t_3)$.

largest spectral radii among $\mathcal{T}(n,k)$, You et al. [12] proposed the following Conjecture 1.

**Conjecture 1** *[12] Let $\mathcal{T} \in \mathcal{T}(n,k) \backslash \mathcal{T}^{(8)}(n,k)$, where $k \geq 3$. For $m = \frac{n-1}{k-1} \geq 10$, we have*

$$\rho_\alpha\left(\mathcal{S}_{m+1}^k\right) > \rho_\alpha\left(\mathcal{S}_{1,m-2}^k\right) > \rho_\alpha\left(\mathcal{S}_{2,m-3}^k\right) > \rho_\alpha(\mathcal{T}(1,1,m-3)) >$$
$$\rho_\alpha\left(\mathcal{S}_{1,m-4,1}^k\right) > \rho_\alpha\left(\mathcal{S}_{m-3,0,1}^k\right) > \rho_\alpha\left(\mathcal{S}_{3,m-4}^k\right) > \rho_\alpha(\mathcal{T}(1,2,m-4)) > \rho_\alpha(\mathcal{T}). \quad (1)$$

Conjecture 1 was partially solved for $1 - \frac{1}{m-2} \leq \alpha < 1$ and $m \geq 10$ by Wang et al. [13]. However, for $0 \leq \alpha < 1 - \frac{1}{m-2}$ and $m \geq 10$, whether Conjecture 1 is correct or not remains a problem to be further solved. The goal of this article is to completely prove that Conjecture 1 is right for $0 \leq \alpha < 1$. This article is organized as follows. In Section 2, some lemmas which are necessary for the subsequent proofs are introduced. In Section 3, to obtain our results, we introduce an effective and new $\rho_\alpha$-normal labeling method for the $\alpha$-spectral radius of the $k$-uniform hypergraphs. For a connected $k$-uniform hypergraph, we generalize the definition of the weighted incidence matrix for its spectral radius developed by Lu and Man [14] to its $\alpha$-spectral radius. The relationship between the weighted incidence matrix and the $\alpha$-spectral radius of a connected $k$-uniform hypergraph is characterized. In Section 4, by using the new $\rho_\alpha$-normal labeling method proposed in this article and the majorization theorem derived by You et al. [12], we prove that Conjecture 1 is correct for $0 \leq \alpha < 1$ and $m \geq 13$.

## 2. PRELIMINARIES

In this section, some definitions and necessary lemmas are introduced.

**Definition 2.1** *[15, 16] Let $\mathcal{A} = (a_{i_1 i_2 \cdots i_k})$ be a nonnegative tensor of order $k$ and dimension $s$. For any nonempty proper index subset $I \subset [s] = \{1, 2, \cdots, s\}$, if there is at least an entry $a_{i_1 i_2 \cdots i_k} > 0$, where $i_1 \in I$ and at least an $i_j \in [s] \setminus I$ for $j = 2, 3, \ldots, k$, then $\mathcal{A}$ is called a nonnegative weakly irreducible tensor.*

Let $\mathbb{R}_+^n = \{\boldsymbol{x} \in \mathbb{R}^n \mid x_i \geq 0, \forall i \in [n]\}$ and $\mathbb{R}_{++}^n = \{\boldsymbol{x} \in \mathbb{R}^n \mid x_i > 0, \forall i \in [n]\}$, where $[n] = \{1, 2, \cdots, n\}$.

**Lemma 2.1** *[15, 17] (The Perron–Frobenius theorem for nonnegative tensors) Let $\mathcal{A}$ be a nonnegative tensor of order $k$ and dimension $n$, where $k \geq 2$. Then we have the following statements.*

*(i). $\rho(\mathcal{A})$ is an eigenvalue of $\mathcal{A}$ with a nonnegative eigenvector $\boldsymbol{x} \in \mathbb{R}_+^n$ corresponding to it.*

*(ii). If $\mathcal{A}$ is weakly irreducible, then $\rho(\mathcal{A})$ is the only eigenvalue of $\mathcal{A}$ with a positive eigenvector $\boldsymbol{x} \in \mathbb{R}_{++}^n$, up to a positive scaling coefficient.*

**Lemma 2.2** *[18] A $k$-uniform hypergraph $\mathcal{G}$ is connected if and only if $\mathcal{A}_\alpha(\mathcal{G})$ is weakly irreducible.*

From Lemmas 2.1 and 2.2, we know that if $\mathcal{G}$ is a connected $k$-uniform hypergraph, then there exists a unique $\boldsymbol{x} \in \mathbb{R}_{++}^n$ with $\|\boldsymbol{x}\|_k^k = 1$ corresponding to $\rho_\alpha(\mathcal{G})$. This unique vector $\boldsymbol{x}$ is called the $\alpha$-Perron vector of $\mathcal{G}$.

**Lemma 2.3** *[19] Let $\mathcal{A}$ be a nonnegative symmetric tensor of order $k$ and dimension $n$. Then we have*

$$\rho(\mathcal{A}) = \max\left\{ \boldsymbol{x}^T(\mathcal{A}\boldsymbol{x}) \mid \boldsymbol{x} \in \mathbb{R}_+^n, \|\boldsymbol{x}\|_k^k = 1 \right\}.$$

*Furthermore, $\boldsymbol{x} \in \mathbb{R}_+^n$ with $\|\boldsymbol{x}\|_k^k = 1$ is an optimal solution of the above optimization problem if and only if it is an eigenvector of $\mathcal{A}$ corresponding to the eigenvalue $\rho(\mathcal{A})$.*

From Lemma 2.3, we obtain that $\rho(\mathcal{A}_\alpha)$ can also be expressed as follows.

$$\rho(\mathcal{A}_\alpha) = \max\left\{ \frac{\boldsymbol{x}^{\mathrm{T}}(\mathcal{A}_\alpha \boldsymbol{x})}{\|\boldsymbol{x}\|_k^k}, \boldsymbol{x} \in \mathbb{R}_+^n, \boldsymbol{x} \neq 0 \right\}. \tag{2}$$

Let $\mathcal{G} = (V(\mathcal{G}), E(\mathcal{G}))$ be a $k$-uniform hypergraph. Let $\boldsymbol{x} = (x_1, x_2, \cdots, x_n)$ be a vector of dimension $n$ and $V \subseteq V(\mathcal{G})$. We write $x^V = \prod_{i \in V} x_i$ for short. When there is no scope for ambiguity, we write $d_v$ instead of $d_{\mathcal{G}}(v)$. Then we have

$$\boldsymbol{x}^{\mathrm{T}}(\mathcal{A}_\alpha(\mathcal{G})\boldsymbol{x}) = \alpha \sum_{v \in V(\mathcal{G})} d_v x_v^k + (1-\alpha) k \sum_{e: v \in e} x^e. \tag{3}$$

Let $\boldsymbol{x}$ be the $\alpha$-Perron vector of $\mathcal{A}_\alpha(\mathcal{G})$. By Lemma 2.3, the eigenequation of $\mathcal{G}$ at a vertex $v$ of $\mathcal{G}$ for $\rho_\alpha(\mathcal{G})$ is:

$$\rho_\alpha(\mathcal{G}) x_v^{k-1} = (\mathcal{A}_\alpha(\mathcal{G})\boldsymbol{x})_v = \alpha d_v x_v^{k-1} + (1-\alpha) \sum_{e: v \in e} x^{e \setminus \{v\}}, \text{ for each } v \in V(\mathcal{G}). \tag{4}$$

Let $\pi = \{d_0, d_1, \ldots, d_{n-1}\}$ be a non-increasing degree sequence of a $k$-uniform supertree $\mathcal{T}$ with $n$ vertices, where $d_0 \geq d_1 \geq \ldots \geq d_{n-1}$. Let

$$\mathcal{T}_\pi = \{\mathcal{T} \mid \mathcal{T} \text{ is a } k\text{-uniform supertree on n vertices with } \pi \text{ as its degree sequence}\}.$$

For a $k$-uniform supertree $\mathcal{T} = (V(\mathcal{T}), E(\mathcal{T}))$, we take a vertex $v_0$ of $V(\mathcal{T})$ as the root vertex of $\mathcal{T}$. Let $v \in V(\mathcal{T})$. The minimum length of a path connecting $v$ and $v_0$ is defined as the height of $v$ and is denoted by $h(v)$. We also say that $v$ is in the layer $h(v)$ of $\mathcal{T}$. Let $h = \max_{v \in V(\mathcal{T})} \{h(v)\}$. Xiao et al. [20] generalized the breadth-first-search ordering (BFS-ordering for short) on trees to $k$-uniform supertrees, which is shown in Definition 2.2.

**Definition 2.2** *[20] Let $\mathcal{T} = (V(\mathcal{T}), E(\mathcal{T}))$ be a $k$-uniform supertree with root $v_0$, where $k \geq 2$. An ordering $\prec$ of the vertices of $\mathcal{T}$ is called a BFS-ordering if all the following (i)–(iv) hold for all vertices of $\mathcal{T}$:*

(i). $u \prec v$ implies $h(u) \leq h(v)$.

(ii). $u \prec v$ implies $d_u \geq d_v$.

(iii). If $\{u, u_1\} \subset e_1 \in E(\mathcal{T})$ and $\{v, v_1\} \subset e_2 \in E(\mathcal{T})$ such that $u \prec v$, $h(u) = h(u_1) + 1$, and $h(v) = h(v_1) + 1$, then $u_1 \prec v_1$.

(iv). Suppose $u_1 \prec u_2 \prec \cdots \prec u_r$ for every edge $e = \{u_1, u_2, \ldots, u_k\} \in E(\mathcal{T})$, then there exist no vertices $v \in V(\mathcal{T}) \backslash e$ such that $u_i \prec v \prec u_{i+1}$, where $2 \leq i \leq k-1$.

**Lemma 2.4** [12] Let $\pi = \{d_0, d_1, \ldots, d_{n-1}\}$ be a degree sequence of a k-uniform supertree having $n$ vertices with $d_0 \geq d_1 \geq \cdots \geq d_{n-1}$. Let $\mathcal{T}$ be a k-uniform supertree with the largest $\alpha$-spectral radius in $\mathcal{T}_\pi$. Then $\mathcal{T}$ has a BFS-ordering.

Xiao et al. [20] constructed a BFS-supertree $\mathcal{T}^* \in \mathcal{T}_\pi$ with a degree sequence when they studied the supertree with the largest spectral radius among the k-uniform supertrees with a given degree sequence. Later, You et al. [12] improved the construction of a BFS-supertree $\mathcal{T}^*$ with a degree sequence $\pi$. You et al. [12] obtained that for any given degree sequence $\pi$, there exists a k-uniform supertree $\mathcal{T}^*$ with the degree sequence $\pi$ and having a BFS-ordering. For more details of the construction for $\mathcal{T}^*$, one can refer to [12, 20].

**Lemma 2.5** [20] Let $\mathcal{T} \in \mathcal{T}_\pi$ with a BFS-ordering. Then $\mathcal{T} \cong \mathcal{T}^*$.

**Lemma 2.6** [12] Let $\mathcal{T} \in \mathcal{T}_\pi$. Then $\mathcal{T}$ has the largest $\alpha$-spectral radius in $\mathcal{T}_\pi$ if and only if $\mathcal{T} \cong \mathcal{T}^*$.

**Definition 2.3** [12] Let $\pi = (d_0, d_1, \ldots, d_{n-1})$ and $\pi' = (d'_0, d'_1, \ldots, d'_{n-1})$ be two non-increasing degree sequences. We say that $\pi$ is majorized by $\pi'$, denoted by $\pi \triangleleft \pi'$, if $\sum_{i=0}^{n-1} d_i = \sum_{i=0}^{n-1} d'_i$ and $\sum_{i=0}^{j} d_i \leq \sum_{i=0}^{j} d'_i$ for all $j = 0, 1, 2, \ldots, n-2$.

**Lemma 2.7** [12] Let $\pi = (d_0, d_1, \ldots, d_{n-1})$ and $\pi' = (d'_0, d'_1, \ldots, d'_{n-1})$ be two different non-increasing degree sequences. Suppose that $\mathcal{T} = (V(\mathcal{T}), E(\mathcal{T}))$ and $\mathcal{T}' = (V(\mathcal{T}'), E(\mathcal{T}'))$ attain the largest $\alpha$-spectral radii in $\mathcal{T}_\pi$ and $\mathcal{T}_{\pi'}$, respectively. If $\pi \triangleleft \pi'$, then $\rho_\alpha(\mathcal{T}') > \rho_\alpha(\mathcal{T})$.

**Lemma 2.8** [12] Let $k \geq 3$ and $m = \frac{n-1}{k-1} \geq 7$ be positive integers. Then

(i). $\rho_\alpha\left(\mathcal{S}_{m+1}^k\right) > \rho_\alpha\left(\mathcal{S}_{1,m-2}^k\right) > \rho_\alpha\left(\mathcal{S}_{2,m-3}^k\right) > \cdots > \rho_\alpha\left(\mathcal{S}_{\lfloor \frac{m-1}{2} \rfloor, \lceil \frac{m-1}{2} \rceil}^k\right)$.

(ii). For any supertree $\mathcal{T} \in \mathcal{T}(n,k) \backslash \{\mathcal{S}_{m+1}^k, \mathcal{S}_{1,m-2}^k, \mathcal{S}_{2,m-3}^k\}$, $\rho_\alpha\left(\mathcal{S}_{2,m-3}^k\right) > \rho_\alpha(\mathcal{T})$.

# 3. A NEW $\rho_\alpha$-NORMAL LABELING METHOD FOR THE $\alpha$-SPECTRAL RADIUS OF $k$-UNIFORM HYPERGRAPHS

In this section, to solve Conjecture 1, we will propose a useful $\rho_\alpha$-normal labeling method for the $\alpha$-spectral radius of the $k$-uniform hypergraphs, which generalizes the $\alpha$-normal labeling method developed by Lu and Man [14] for the spectral radius of the $k$-uniform hypergraphs. The definitions of $\rho_\alpha$-normal, $\rho_\alpha$-subnormal and $\rho_\alpha$-supernormal for the $\alpha$-spectral radius of the $k$-uniform hypergraphs are introduced, which are shown in Definitions 3.1–3.3, respectively. Then, we characterize the relationship between the $\rho_\alpha$-normal labeling and the $\alpha$-spectral radius of $k$-uniform hypergraphs, which are shown in Lemmas 3.1–3.3.

**Definition 3.1** *Let $k \geq 2$ and $0 \leq \alpha < 1$. A connected $k$-uniform hypergraph $\mathcal{G}$ is called $\rho_\alpha$-normal if there exists a weighted incidence matrix $\boldsymbol{B}$ satisfying*

*(i). $\sum_{e:v\in e} \big(B(v,e) + \alpha\big) = \rho_\alpha$, for any $v \in V(\mathcal{G})$.*

*(ii). $\prod_{v:v\in e} B(v,e) = (1-\alpha)^k$, for any $e \in E(\mathcal{G})$.*

*Moreover, the incidence matrix $\boldsymbol{B}$ is called consistent if for any cycle $v_0 e_1 v_1 \ldots v_l$ ($v_0 = v_l$) of $\mathcal{G}$, we have $\prod_{i=1}^{l} \frac{B(v_i,e_i)}{B(v_{i-1},e_i)} = 1$. In this case, we call $\mathcal{G}$ consistently $\rho_\alpha$-normal.*

**Remark 3.1** *For any supertree $\mathcal{T}$, since $\mathcal{T}$ does not contain cycles, $\mathcal{T}$ satisfies the consistent condition naturally.*

**Lemma 3.1** *Let $\mathcal{G}$ be a connected $k$-uniform hypergraph, where $k \geq 2$. The $\alpha$-spectral radius of $\mathcal{G}$ is $\rho_\alpha$ if and only if $\mathcal{G}$ is consistently $\rho_\alpha$-normal, where $0 \leq \alpha < 1$.*

**Proof.** Let $V(\mathcal{G}) = \{v_1, v_2, \cdots, v_n\}$.

(1). The proof of necessity.

We suppose that the $\alpha$-spectral radius of $\mathcal{G}$ is $\rho_\alpha$. We will prove that $\mathcal{G}$ is consistently $\rho_\alpha$-normal. Let $\boldsymbol{x} = (x_1, x_2, \cdots, x_n)$ be the $\alpha$-Perror eigenvector of the $\alpha$-spectral radius of $\mathcal{G}$. We define the weighted incidence matrix $\boldsymbol{B}$ as follows. Let

$$B(v,e) = \begin{cases} \frac{(1-\alpha)x^e}{x_v^k}, & \text{if } v \in e, \\ 0, & \text{otherwise.} \end{cases}$$

(1.1). For any $v \in V(\mathcal{G})$, we have

$$\sum_{e:v\in e} \big(B(v,e) + \alpha\big) = \sum_{e:v\in e} \left(\frac{(1-\alpha)x^e}{x_v^k} + \alpha\right) = \frac{\alpha d_v x_v^k + (1-\alpha)\sum_{e:v\in e} x^e}{x_v^k}. \tag{5}$$

By the eigenequation (4) of $\mathcal{G}$ at $v$, we get

$$\rho_\alpha x_v^k = \alpha d_v x_v^k + (1-\alpha) \sum_{e:v\in e} x^e. \tag{6}$$

Therefore, by substituting (6) into (5), we get $\sum_{e:v\in e} \big(B(v,e) + \alpha\big) = \rho_\alpha$. Namely, we have Definition 3.1 (i).

(1.2). For any $e \in E(\mathcal{G})$, we get

$$\prod_{v:v\in e} B(v,e) = \prod_{v:v\in e} \frac{(1-\alpha)x^e}{x_v^k} = (1-\alpha)^k \cdot \frac{(x^e)^k}{\prod_{v:v\in e} x_v^k} = (1-\alpha)^k, \tag{7}$$

where the third equality in (7) holds since $\prod_{v:v\in e} x_v^k = (x^e)^k$. By (7), we have Definition 3.1 (ii).

Next, we prove that $\boldsymbol{B}$ is consistent. For any cycle $v_0 e_1 v_1 \ldots v_l$ ($v_l = v_0$) of $\mathcal{G}$, we obtain

$$\prod_{i=1}^{l} \frac{B(v_i, e_i)}{B(v_{i-1}, e_i)} = \prod_{i=1}^{l} \frac{\frac{(1-\alpha)x^{e_i}}{x_{v_i}^k}}{\frac{(1-\alpha)x^{e_i}}{x_{v_{i-1}}^k}} = \prod_{i=1}^{l} \frac{x_{v_{i-1}}^k}{x_{v_i}^k} = \frac{x_{v_0}^k}{x_{v_l}^k} = 1. \tag{8}$$

By (8), we get that $\mathcal{G}$ is consistently $\rho_\alpha$-normal.

(2). The proof of sufficiency.

Suppose that $\mathcal{G}$ is consistently $\rho_\alpha$-normal. We will prove that the $\alpha$-spectral radius of $\mathcal{G}$ is $\rho_\alpha$. Let $\boldsymbol{x} = (x_1, \ldots, x_n)$ be an arbitrary nonzero vector in $\mathbb{R}_+^n$.

For any $e \in E(\mathcal{G})$, if $\prod_{v:v\in e} B(v,e) = (1-\alpha)^k$, then we have

$$(1-\alpha) \sum_{e\in E(\mathcal{G})} \frac{k}{1-\alpha} \prod_{v:v\in e} \left((B(v,e))^{\frac{1}{k}} x_v\right) = (1-\alpha) \sum_{e\in E(\mathcal{G})} kx^e. \tag{9}$$

By the Arithmetic Mean–Geometry Mean inequality, we get

$$\sum_{e\in E(\mathcal{G})} k \prod_{v:v\in e} \left(B(v,e)^{\frac{1}{k}} x_v\right) \le \frac{\sum_{e\in E(\mathcal{G})} \sum_{v:v\in e} k B(v,e) x_v^k}{k}. \tag{10}$$

Obviously, we have

$$\alpha \sum_{v\in V(\mathcal{G})} d_v x_v^k = \sum_{v\in V(\mathcal{G})} \sum_{e:v\in e} \alpha x_v^k. \tag{11}$$

By (3), (9)–(11) and Condition (i) in Definition 3.1, we have

$$\boldsymbol{x}^T (\boldsymbol{\mathcal{A}}_\alpha(\mathcal{G})\boldsymbol{x}) = \alpha \sum_{v\in V(\mathcal{G})} d_v x_v^k + (1-\alpha) \sum_{e\in E(\mathcal{G})} kx^e$$

$$\leq \sum_{v \in V(\mathcal{G})} \sum_{e:v \in e} \big(\alpha + B(v,e)\big) x_v^k = \rho_\alpha \sum_{v \in V(\mathcal{G})} x_v^k = \rho_\alpha \parallel \boldsymbol{x} \parallel_k^k. \tag{12}$$

Therefore, by (12) and the arbitrariness of $\boldsymbol{x}$, we obtain $\rho_\alpha(\mathcal{G}) \leq \rho_\alpha$, with the equality if and only if $\mathcal{G}$ is $\rho_\alpha$-normal and the equality in (10) holds. Namely, there is a nonzero solution $\{x_i\}$ for the system of the following homogeneous linear equations:

$$B(v_{i_1},e)^{\frac{1}{k}} x_{v_{i_1}} = B(v_{i_2},e)^{\frac{1}{k}} x_{v_{i_2}} = \ldots = B(v_{i_k},e)^{\frac{1}{k}} x_{v_{i_k}}, \forall e = \{v_{i_1}, \ldots, v_{i_k}\} \in E(\mathcal{G}). \tag{13}$$

Let $v_0$ be an arbitrary vertex in $V(\mathcal{G})$. For any $u \in V(\mathcal{G})$, since $\mathcal{G}$ is connected, there exists a path $v_0 e_1 v_1 e_2 v_2 \ldots v_l$ ($v_l = u$) connecting $v_0$ and $u$. Let $x_{v_0}^* = 1$. For $u \in V(\mathcal{G})$, we define $x_u^* = \left( \prod_{i=1}^{l} \frac{B(v_{i-1},e_i)}{B(v_i,e_i)} \right)^{\frac{1}{k}}$. The consistent condition guarantees that $x_u^*$ is independent of the choice of the path. We can check that $(x_1^*, x_2^*, \cdots, x_n^*)$ is a solution of (13). Thus, we have $\rho_\alpha(\mathcal{G}) = \rho_\alpha$. $\square$

**Definition 3.2** *Let $k \geq 2$ and $0 \leq \alpha < 1$. A connected $k$-uniform hypergraph $\mathcal{G}$ is called $\rho_\alpha$-subnormal if there exists a weighted incidence matrix $\boldsymbol{B}$ satisfying*

*(i).* $\sum_{e:v \in e} \big(B(v,e) + \alpha\big) \leq \rho_\alpha$, *for any* $v \in V(\mathcal{G})$.
*(ii).* $\prod_{v:v \in e} B(v,e) \geq (1-\alpha)^k$, *for any* $e \in E(\mathcal{G})$.

*Moreover, $\mathcal{G}$ is called strictly $\rho_\alpha$-subnormal if it is $\rho_\alpha$-subnormal but not $\rho_\alpha$-normal.*

**Lemma 3.2** *Let $\mathcal{G}$ be a connected $k$-uniform hypergraph, where $k \geq 2$. If $\mathcal{G}$ is $\rho_\alpha$-subnormal, then $\rho_\alpha(\mathcal{G}) \leq \rho_\alpha$, where $0 \leq \alpha < 1$. Moreover, if $\mathcal{G}$ is strictly $\rho_\alpha$-subnormal, then $\rho_\alpha(\mathcal{G}) < \rho_\alpha$.*

**Proof**. Let $\boldsymbol{x} = (x_1, \ldots, x_n)$ be an arbitrary nonzero vector in $\mathbb{R}_+^n$. For any $e \in E(\mathcal{G})$, if $\prod_{v:v \in e} B(v,e) \geq (1-\alpha)^k$, then we have

$$(1-\alpha) \sum_{e \in E(\mathcal{G})} \frac{k}{1-\alpha} \prod_{v:v \in e} \left( (B(v,e))^{\frac{1}{k}} x_v \right) \geq (1-\alpha) \sum_{e \in E(\mathcal{G})} kx^e. \tag{14}$$

By (3), (10), (11), (14), and Condition (i) in Definition 3.2, we have

$$\boldsymbol{x}^T (\boldsymbol{\mathcal{A}}_\alpha(\mathcal{G}))\boldsymbol{x} = \alpha \sum_{v \in V(\mathcal{G})} d_v x_v^k + (1-\alpha) \sum_{e \in E(\mathcal{G})} kx^e$$

$$\leq \sum_{v \in V(\mathcal{G})} \sum_{e:v \in e} \big(\alpha + B(v,e)\big) x_v^k \leq \rho_\alpha \sum_{v \in V(\mathcal{G})} x_v^k = \rho_\alpha \parallel \boldsymbol{x} \parallel_k^k. \tag{15}$$

Therefore, by (15) and the arbitrariness of $\boldsymbol{x}$, we obtain $\rho_\alpha(\mathcal{G}) \leq \rho_\alpha$. If $\mathcal{G}$ is strictly $\rho_\alpha$-subnormal, then the inequality in (14) or the second inequality in (15) holds. Thus, we get $\rho_\alpha(\mathcal{G}) < \rho_\alpha$. $\square$

**Definition 3.3** *Let $k \geq 2$ and $0 \leq \alpha < 1$. A connected $k$-uniform hypergraph $\mathcal{G}$ is called $\rho_\alpha$-supernormal if there exists a weighted incidence matrix $\boldsymbol{B}$ satisfying*

*(i).* $\sum_{e:v\in e} \big(B(v,e) + \alpha\big) \geq \rho_\alpha$, *for any $v \in V(\mathcal{G})$.*

*(ii).* $\prod_{v:v\in e} B(v,e) \leq (1-\alpha)^k$, *for any $e \in E(\mathcal{G})$.*

*Moreover, $\mathcal{G}$ is called strictly $\rho_\alpha$-supernormal if it is $\rho_\alpha$-supernormal but not $\rho_\alpha$-normal.*

**Lemma 3.3** *Let $\mathcal{G}$ be a connected $k$-uniform hypergraph, where $k \geq 2$. If $G$ is consistently $\rho_\alpha$-supernormal, then $\rho_\alpha(\mathcal{G}) \geq \rho_\alpha$, where $0 \leq \alpha < 1$. Moreover, if $\mathcal{G}$ is strictly consistently $\rho_\alpha$-supernormal, then $\rho_\alpha(\mathcal{G}) > \rho_\alpha$.*

**Proof**. From the consistent condition of $\mathcal{G}$ and the proof of sufficiency of Lemma 3.1, there exists an eigenvector $\boldsymbol{x} = (x_1, x_2, \cdots, x_n) \in \mathbb{R}^n_+$ such that (13) holds. We have

$$\sum_{e\in E(\mathcal{G})} \prod_{v:v\in e} \left(B(v,e)^{\frac{1}{k}} x_v\right) = \sum_{e\in E(\mathcal{G})} \frac{\sum_{v:v\in e} B(v,e) x_v^k}{k}. \tag{16}$$

For any $e \in E(\mathcal{G})$, if $\prod_{v:v\in e} B(v,e) \leq (1-\alpha)^k$, then

$$(1-\alpha) \sum_{e\in E(\mathcal{G})} \frac{k}{1-\alpha} \prod_{v:v\in e} \left((B(v,e))^{\frac{1}{k}} x_v\right) \leq (1-\alpha) \sum_{e\in E(\mathcal{G})} kx^e. \tag{17}$$

By (3), (11), (16), (17), and Condition (i) in Definition 3.3, we get

$$\boldsymbol{x}^T(\boldsymbol{\mathcal{A}}_\alpha(\mathcal{G})\boldsymbol{x}) = \alpha \sum_{v\in V(\mathcal{G})} d_v x_v^k + (1-\alpha) \sum_{e\in E(\mathcal{G})} kx^e$$
$$\geq \sum_{v\in V(\mathcal{G})} \sum_{e:v\in e} \big(\alpha + B(v,e)\big) x_v^k \geq \rho_\alpha \sum_{v\in V(\mathcal{G})} x_v^k = \rho_\alpha \parallel \boldsymbol{x} \parallel_k^k. \tag{18}$$

Therefore, by (18), we obtain $\rho_\alpha(\mathcal{G}) \geq \frac{\boldsymbol{x}^T(\boldsymbol{\mathcal{A}}_\alpha \boldsymbol{x})}{\|\boldsymbol{x}\|_k^k} \geq \rho_\alpha$. If $\mathcal{G}$ is strictly $\rho_\alpha$-supernormal, then the inequality in (17) or the second inequality in (18) holds. Thus, we get $\rho_\alpha(\mathcal{G}) > \rho_\alpha$. □

## 4. A COMPLETE SOLUTION OF THE CONJECURE ON THE $k$-UNIFORM SUPERTREES WITH THE EIGHT LARGEST $\alpha$-SPECTRAL RADII

In this section, we will apply the new $\rho_\alpha$-normal labeling method proposed in Section 3 to compare the $\alpha$-spectral radii of the hypergraphs in Conjecture 1. We introduce Lemmas 4.1–4.3 firstly. Then, by Lemmas 4.1–4.3 and the majorization theorem derived by You et al. [12], we obtain that Conjecture 1 is right for $0 \leq \alpha < 1$, $k \geq 3$, and $m = \frac{n-1}{k-1} \geq 13$, which is shown in Theorem 4.1.

**Lemma 4.1** *We have $\rho_\alpha\big(\mathcal{T}(1,1,m-3)\big) > \rho_\alpha\big(\mathcal{S}_{1,m-4,1}^k\big) > \rho_\alpha\big(\mathcal{S}_{m-3,0,1}^k\big)$, where $k \geq 3$ and $m = \frac{n-1}{k-1} \geq 7$.*

**Proof.** Let $\rho_\alpha^* = \rho_\alpha\big(\mathcal{S}_{1,m-4,1}^k\big)$. Let $\boldsymbol{x} = (x_1, \ldots, x_n) \in \mathbb{R}_{++}^n$ be the $\alpha$-Perron vector of $\rho_\alpha^*$. We suppose that $\mathcal{S}_{1,m-4,1}^k$ is consistently $\rho_\alpha^*$-normal. In $\mathcal{S}_{s_1,s_2,s_3}^k$ (as shown in Fig. 1(b)), let $s_1 = s_3 = 1$ and $s_2 = m-4$ and we get $\mathcal{S}_{1,m-4,1}^k$. Let $v_1$, $v_2$ and $v_3$ be three vertices of degree 1 in $\mathcal{S}_{1,m-4,1}^k$, where $v_1$, $v_2$ and $v_3$ are respectively incident with $e_1$, $e_3$ and $e_4$ of $\mathcal{S}_{1,m-4,1}^k$. The vertices $v_1$, $v_2$ and $v_3$ and the edges $e_1$, $e_3$ and $e_4$ of $\mathcal{S}_{1,m-4,1}^k$ are shown in Fig. 1(b). By using the eigenequations (4) of $\mathcal{S}_{1,m-4,1}^k$ at $v_1$, $v_2$, $v_3$, $u_1$, and $u_2$, and bearing the symmetry of the entries in $\boldsymbol{x}$ in mind, we get respectively

$$\rho_\alpha^* x_{v_1}^{k-1} = \alpha x_{v_1}^{k-1} + (1-\alpha) x_{v_1}^{k-3} x_{u_1} x_{u_2}, \tag{19}$$

$$\rho_\alpha^* x_{v_2}^{k-1} = \alpha x_{v_2}^{k-1} + (1-\alpha) x_{v_2}^{k-2} x_{u_1}, \tag{20}$$

$$\rho_\alpha^* x_{v_3}^{k-1} = \alpha x_{v_3}^{k-1} + (1-\alpha) x_{v_3}^{k-2} x_{u_2}, \tag{21}$$

$$\rho_\alpha^* x_{u_1}^{k-1} = 2\alpha x_{u_1}^{k-1} + (1-\alpha) x_{v_2}^{k-1} + (1-\alpha) x_{v_1}^{k-2} x_{u_2}, \tag{22}$$

$$\rho_\alpha^* x_{u_2}^{k-1} = (m-2)\alpha x_{u_2}^{k-1} + (m-4)(1-\alpha) x_{v_3}^{k-1} + 2(1-\alpha) x_{v_1}^{k-2} x_{u_1}. \tag{23}$$

From (19), we have $\rho_\alpha^* > \alpha$ when $\boldsymbol{x} \in \mathbb{R}_{++}^n$ and $0 \leq \alpha < 1$. Let $A_0 = \frac{1-\alpha}{\rho_\alpha^* - \alpha}$. Therefore, we have $A_0 > 0$. By (19)–(21), we get respectively

$$x_{v_1} = \sqrt{A_0 x_{u_1} x_{u_2}}, \quad x_{v_2} = A_0 x_{u_1}, \quad x_{v_3} = A_0 x_{u_2}. \tag{24}$$

For simplicity, let

$$A_1 = A_0^k, \quad B_0 = (1-\alpha) A_0^{k-1}, \quad B_1 = \rho_\alpha^* - 2\alpha - B_0, \quad C_1 = \frac{(\rho_\alpha^* - \alpha) B_0}{B_1}. \tag{25}$$

Obviously, since $0 \leq \alpha < 1$ and $A_0 > 0$, we get $B_0 > 0$.

Substituting (24) into (22), we get

$$\rho_\alpha^* x_{u_1}^{k-1} = 2\alpha x_{u_1}^{k-1} + (1-\alpha)(A_0 x_{u_1})^{k-1} + (1-\alpha)(A_0 x_{u_1} x_{u_2})^{\frac{k-2}{2}} x_{u_2}. \tag{26}$$

From (26), we obtain

$$\rho_\alpha^* > 2\alpha + B_0, \quad (\text{since } 0 \leq \alpha < 1 \text{ and } A_0 > 0), \tag{27}$$

$$x_{u_1} = \frac{(1-\alpha)}{B_1^{\frac{2}{k}} (\rho_\alpha^* - \alpha)^{1-\frac{2}{k}}} x_{u_2}. \tag{28}$$

Substituting (24) and (28) into (23), we get

$$\rho_\alpha^* - (m-2)\alpha = (m-4)B_0 + 2C_1. \tag{29}$$

By (25) and (27), we get $B_1 > 0$. Since $\rho_\alpha^* > \alpha$ and $B_0, B_1 > 0$, by (25), we get $C_1 > 0$. Since $B_0, C_1 > 0$ and $0 \leq \alpha < 1$, by (29), when $m \geq 5$, we obtain

$$\rho_\alpha^* - 2\alpha - C_1 = (m-4)\alpha + (m-4)B_0 + C_1 > 0. \tag{30}$$

(1). The proof of $\rho_\alpha\big(\mathcal{T}(1,1,m-3)\big) > \rho_\alpha(\mathcal{S}_{1,m-4,1}^k)$ for $m = \frac{n-1}{k-1} \geq 7$

In $\mathcal{T}(t_1, t_2, t_3)$ (as shown in Fig. 1(c)), let $t_1 = t_2 = 1$ and $t_3 = m-3$. Namely, we get $\mathcal{T}(1,1,m-3)$. We construct a weighted incidence matrix $\boldsymbol{B}$ for $\mathcal{T}(1,1,m-3)$ as follows. Let $B(v, e_i) = \rho_\alpha^* - \alpha$, where $v$ is an arbitrary vertex of degree 1 in $V(\mathcal{T}(1,1,m-3))$ and $e_i$ $(1 \leq i \leq m)$ is the edge incident with $v$. Let

$$B(u_1, e_1) = B(u_2, e_1) = B_1, \quad B(u_3, e_1) = \frac{(1-\alpha)^3 A_0^{k-3}}{B_1^2}.$$

Furthermore, let $B(u_1, e_2) = B(u_2, e_3) = B(u_3, e_i) = B_0$, where $4 \leq i \leq m$. Since $0 \leq \alpha < 1$ and $A_0, B_0, B_1 > 0$, we get $B(v, e) > 0$ for any vertex $v$ and any edge $e$ incident with $v$ in $\mathcal{T}(1,1,m-3)$.

It is easy to verify that $\prod_{v:v \in e_i} B(v, e_i) = (1-\alpha)^k$ for any edge $e_i$ $(1 \leq i \leq m)$ in $E\big(\mathcal{T}(1,1,m-3)\big)$. We also can verify that $\sum_{e:v \in e}(B(v,e) + \alpha) = \rho_\alpha^*$, where $v$ is an arbitrary vertex of degree 1 in $\mathcal{T}(1,1,m-3)$ or $v = u_1, u_2$.

Next, we compare $\sum_{e:u_3 \in e}\big(B(u_3, e) + \alpha\big)$ with $\rho_\alpha^*$. By (29), we have

$$\rho_\alpha^* - \sum_{e:u_3 \in e}\big(B(u_3, e) + \alpha\big)$$
$$= \rho_\alpha^* - \big(B(u_3, e_1) + (m-3)B(u_3, e_4) + (m-2)\alpha\big)$$
$$= 2C_1 - B_0 - \frac{(1-\alpha)^3 A_0^{k-3}}{B_1^2}$$
$$= -(1-\alpha)^3 A_0^{k-3}\left(\frac{1}{B_1} - \frac{1}{\rho_\alpha^* - \alpha}\right)^2. \tag{31}$$

Since $B_0 > 0$ and $0 \leq \alpha < 1$, we have $B_1 = \rho_\alpha^* - 2\alpha - B_0 \neq \rho_\alpha^* - \alpha$. By (31), since $A_0 > 0$, we have $\sum_{e:u_3 \in e}(B(u_3, e) + \alpha) > \rho_\alpha^*$. Since $\mathcal{T}(1,1,m-3)$ is a supertree, it follows from Remark 3.1 that $\boldsymbol{B}$ is consistent. Thus, $\mathcal{T}(1,1,m-3)$ is strictly consistently $\rho_\alpha^*$-supernormal. By Lemma 3.3, we have $\rho_\alpha\big(\mathcal{T}(1,1,m-3)\big) > \rho_\alpha^* = \rho_\alpha\big(\mathcal{S}_{1,m-4,1}^k\big)$.

(2). The proof of $\rho_\alpha\left(\mathcal{S}_{1,m-4,1}^k\right) > \rho_\alpha\left(\mathcal{S}_{m-3,0,1}^k\right)$, where $m = \frac{n-1}{k-1} \geq 7$.

In $\mathcal{S}_{s_1,s_2,s_3}^k$ (as shown in Fig. 1(b)), let $s_1 = m-3$, $s_2 = 0$ and $s_3 = 1$. Namely, we get $\mathcal{S}_{m-3,0,1}^k$. The vertices $u_1$, $u_2$ and $u_3$ and the edges of $\mathcal{S}_{m-3,0,1}^k$ are shown in Fig. 1(b). We construct a weighted incidence matrix $\boldsymbol{B}$ for $\mathcal{S}_{m-3,0,1}^k$ as follows. Let $B(v, e_i) = \rho_\alpha^* - \alpha$, where $v$ is an arbitrary vertex of degree 1 in $V(\mathcal{S}_{m-3,0,1}^k)$ and $e_i$ is the edge incident with $v$. Furthermore, let

$$B(u_1, e_1) = \frac{(1-\alpha)^2 A_0^{k-2}}{\rho_\alpha^* - 2\alpha - C_1}, \qquad B(u_1, e_i) = B_0, \quad \text{for } 3 \leq i \leq m-1;$$

$$B(u_2, e_1) = \rho_\alpha^* - 2\alpha - C_1, \qquad B(u_2, e_2) = C_1;$$

$$B(u_3, e_2) = B_1, \qquad B(u_3, e_m) = B_0.$$

Since $0 \leq \alpha < 1$, $A_0, B_0, B_1, C_1 > 0$, $\rho_\alpha^* - \alpha > 0$, and $\rho_\alpha^* - 2\alpha - C_1 > 0$ (by (30)), we get $B(v, e) > 0$ for any vertex $v$ and any edge $e$ incident with $v$ in $\mathcal{S}_{m-3,0,1}^k$.

It is easy to verify that $\prod_{v:v\in e} B(v,e) = (1-\alpha)^k$ for any edge $e = e_i$ ($1 \leq i \leq m$) in $E(\mathcal{S}_{m-3,0,1}^k)$. We also can verify that $\sum_{e:v\in e}(B(v,e)+\alpha) = \rho_\alpha^*$, where $v$ is an arbitrary vertex of degree 1 in $V\left(\mathcal{S}_{m-3,0,1}^k\right)$ or $v = u_2, u_3$.

Next, we compare $\sum_{e:u_1\in e}\left(B(u_1,e)+\alpha\right)$ with $\rho_\alpha^*$. Since $C_1 > 0$, by (29), we obtain $\rho_\alpha^* \geq 5\alpha + 3B_0$ when $m \geq 7$. Thus, when $m \geq 7$, we get

$$\frac{3}{2}B_0 - C_1 = \frac{B_0}{B_1}\left(\frac{3}{2}B_1 - (\rho_\alpha^* - \alpha)\right) = \frac{B_0}{B_1}\left(\frac{1}{2}\rho_\alpha^* - 2\alpha - \frac{3}{2}B_0\right) \geq 0.$$

From (29) and $C_1 \leq \frac{3}{2}B_0$, when $m \geq 7$, we get

$$\rho_\alpha^* - \sum_{e:u_1\in e}\left(B(u_1,e)+\alpha\right)$$
$$= \rho_\alpha^* - \left(B(u_1,e_1) + (m-3)B(u_1,e_3) + (m-2)\alpha\right)$$
$$= 2C_1 - B_0 - \frac{(1-\alpha)^2 A_0^{k-2}}{\rho_\alpha - 2\alpha - C_1}$$
$$= D_1\left[(\alpha + 2B_0 - C_1)\rho_\alpha^* - 2\alpha^2 - 3B_0\alpha - B_0C_1\right]$$
$$= D_1\left[(m-4)\alpha^2 + (3m-11)B_0\alpha - (m-4)C_1\alpha + 2(m-4)B_0^2 - 2C_1^2 - (m-7)B_0C_1\right]$$
$$\geq D_1\left[(m-4)\alpha^2 + (3m-11)B_0\alpha - \frac{3}{2}(m-4)B_0\alpha + 2(m-4)B_0^2 - \frac{9}{2}B_0^2 - \frac{3}{2}(m-7)B_0^2\right]$$
$$= D_1\left[(m-4)\alpha^2 + (\frac{3}{2}m - 5)B_0\alpha + (\frac{1}{2}m - 2)B_0^2\right], \tag{32}$$

where $D_1 = \frac{C_1}{(\rho_\alpha^* - \alpha)(\rho_\alpha^* - 2\alpha - C_1)}$. Since $C_1 > 0$, $\rho_\alpha^* - \alpha > 0$ and $\rho_\alpha^* - 2\alpha - C_1 > 0$ (by (30)), we have $D_1 > 0$. Therefore, when $m \geq 7$, it follows from (32) that $\rho_\alpha^* > \sum_{e:u_1\in e}\left(B(u_1,e)+\alpha\right)$.

Thus, by Definition 3.2, $\mathcal{S}_{m-3,0,1}^k$ is strictly $\rho_\alpha^*$-subnormal. By Lemma 3.2, $\rho_\alpha\left(\mathcal{S}_{m-3,0,1}^k\right) < \rho_\alpha^* = \rho_\alpha\left(\mathcal{S}_{1,m-4,1}^k\right)$. □

**Lemma 4.2** Let $k \geq 3$ and $m = \frac{n-1}{k-1} \geq 10$. We have $\rho_\alpha\left(\mathcal{S}_{m-3,0,1}^k\right) > \rho_\alpha\left(\mathcal{S}_{3,m-4}^k\right) > \rho_\alpha\left(\mathcal{T}(1,2,m-4)\right)$.

**Proof**. In $\mathcal{S}_{a,b}^k$ (as shown in Fig. 1(a)), let $a = 3$ and $b = m-4$. Namely, we get $\mathcal{S}_{3,m-4}^k$. For simplicity, let $\rho_\alpha^\circ = \rho_\alpha\left(\mathcal{S}_{3,m-4}^k\right)$. Let $\boldsymbol{x} = (x_1,\ldots,x_n) \in \mathbb{R}_{++}^n$ be the $\alpha$-Perron vector of $\rho_\alpha^\circ$. We suppose that $\mathcal{S}_{3,m-4}^k$ is consistently $\rho_\alpha^\circ$-normal. Let $v_1$, $v_2$ and $v_3$ of $\mathcal{S}_{3,m-4}^k$ be three vertices of degree 1, where $v_1$, $v_2$ and $v_3$ are respectively incident with $e_1$, $e_2$ and $e_m$ of $\mathcal{S}_{3,m-4}^k$. The vertices $u_1$, $u_2$, $v_1$, $v_2$, and $v_3$ and the edges $e_1$, $e_2$ and $e_m$ of $\mathcal{S}_{3,m-4}^k$ are shown in Fig. 1(a). By the eigenequations (4) of $\mathcal{S}_{3,m-4}^k$ at $v_1$, $v_2$, $v_3$, $u_1$, and $u_2$ and bearing the symmetry of the entries in $\boldsymbol{x}$ in mind, we get

$$\rho_\alpha^\circ x_{v_1}^{k-1} = \alpha x_{v_1}^{k-1} + (1-\alpha)x_{u_1}x_{u_2}x_{v_1}^{k-3}, \tag{33}$$

$$\rho_\alpha^\circ x_{v_2}^{k-1} = \alpha x_{v_2}^{k-1} + (1-\alpha)x_{u_1}x_{v_2}^{k-2}, \tag{34}$$

$$\rho_\alpha^\circ x_{v_3}^{k-1} = \alpha x_{v_3}^{k-1} + (1-\alpha)x_{u_2}x_{v_3}^{k-2}, \tag{35}$$

$$\rho_\alpha^\circ x_{u_1}^{k-1} = 4\alpha x_{u_1}^{k-1} + 3(1-\alpha)x_{v_2}^{k-1} + (1-\alpha)x_{u_2}x_{v_1}^{k-2}, \tag{36}$$

$$\rho_\alpha^\circ x_{u_2}^{k-1} = (m-3)\alpha x_{u_2}^{k-1} + (m-4)(1-\alpha)x_{v_3}^{k-1} + (1-\alpha)x_{u_1}x_{v_1}^{k-2}. \tag{37}$$

It follows from (33)–(35) that

$$x_{v_1} = \sqrt{\frac{1-\alpha}{\rho_\alpha^\circ - \alpha}x_{u_1}x_{u_2}}, \quad x_{v_2} = \frac{1-\alpha}{\rho_\alpha^\circ - \alpha}x_{u_1}, \quad x_{v_3} = \frac{1-\alpha}{\rho_\alpha^\circ - \alpha}x_{u_2}. \tag{38}$$

Furthermore, from (33), we have $\rho_\alpha^\circ > \alpha$ when $\boldsymbol{x} \in \mathbb{R}_{++}^n$ and $0 \leq \alpha < 1$.

Let

$$A_2 = \frac{(1-\alpha)^k}{(\rho_\alpha^\circ - \alpha)^{k-1}}.$$

Since $\rho_\alpha^\circ > \alpha$ and $0 \leq \alpha < 1$, we have $A_2 > 0$.

Substituting (38) into (36), we get

$$\rho_\alpha^\circ > 4\alpha + 3A_2, \quad (\text{since } \boldsymbol{x} \in \mathbb{R}_{++}^n \text{ and } 0 \leq \alpha < 1), \tag{39}$$

$$x_{u_1} = \frac{1-\alpha}{(\rho_\alpha^\circ - \alpha)^{1-2/k}(\rho_\alpha^\circ - 4\alpha - 3A_2)^{2/k}}x_{u_2}. \tag{40}$$

For simplicity, let

$$B_2 = \rho_\alpha^\circ - 2\alpha - A_2, \quad C_2 = \frac{(\rho_\alpha^\circ - \alpha)A_2}{B_2}.$$

It follows from $0 \leq \alpha < 1$, $A_2 > 0$ and (39) that $B_2 > 0$. Since $A_2, B_2 > 0$ and $\rho_\alpha^\circ > \alpha$, we have $C_2 > 0$.

Substituting (38) and (40) into (37), we obtain

$$\rho_\alpha^\circ - (m-3)\alpha = (m-4)A_2 + \frac{(\rho_\alpha^\circ - \alpha)A_2}{\rho_\alpha^\circ - 4\alpha - 3A_2}. \tag{41}$$

Since $0 \leq \alpha < 1$ and $A_2 > 0$, by (39), we obtain $B_2 = \rho_\alpha^\circ - 2\alpha - A_2 > \rho_\alpha^\circ - 4\alpha - 3A_2 > 0$. Therefore, $C_2 = \frac{(\rho_\alpha^\circ - \alpha)A_2}{B_2} < \frac{(\rho_\alpha^\circ - \alpha)A_2}{\rho_\alpha^\circ - 4\alpha - 3A_2}$. Thus, when $m \geq 10$, from (41), we get

$$\rho_\alpha^\circ - 2\alpha - C_2 = (m-5)\alpha + (m-4)A_2 + \frac{(\rho_\alpha^\circ - \alpha)A_2}{\rho_\alpha^\circ - 4\alpha - 3A_2} - C_2 > 0. \tag{42}$$

(1). The proof of $\rho_\alpha\left(\mathcal{S}_{m-3,0,1}^k\right) > \rho_\alpha\left(\mathcal{S}_{3,m-4}^k\right)$, where $m = \frac{n-1}{k-1} \geq 10$.

We construct a weighted incidence matrix $\boldsymbol{B}$ for $\mathcal{S}_{m-3,0,1}^k$ as follows. The vertices $u_1$, $u_2$, $u_3$, $v_1$, and $v_2$ and the edges $e_1$, $e_2$ and $e_3$ of $\mathcal{S}_{m-3,0,1}^k$ are shown in Fig. 1(b) with $s_1 = m-3$, $s_2 = 0$ and $s_3 = 1$. Let $B(v, e_i) = \rho_\alpha^\circ - \alpha$, where $v$ is an arbitrary vertex of degree 1 in $V(\mathcal{S}_{m-3,0,1}^k)$ and $e_i$ $(1 \leq i \leq m)$ is the edge incident with $v$. Furthermore, let

$$B(u_1, e_1) = \frac{(\rho_\alpha^\circ - \alpha)A_2}{\rho_\alpha^\circ - 2\alpha - C_2}, \qquad B(u_1, e_i) = A_2 \quad \text{for} \quad 3 \leq i \leq m-1;$$

$$B(u_2, e_1) = \rho_\alpha^\circ - 2\alpha - C_2, \qquad B(u_2, e_2) = C_2;$$

$$B(u_3, e_2) = B_2, \qquad B(u_3, e_m) = A_2.$$

Since $A_2, B_2, C_2 > 0$, $\rho_\alpha^\circ - \alpha > 0$, and $\rho_\alpha^\circ - 2\alpha - C_2 > 0$ (by (42)), we have $B(v, e) > 0$ for any vertex $v$ and any edge $e$ incident with $v$ in $\mathcal{S}_{m-3,0,1}^k$.

It can be verified that $\prod_{v:v \in e_i} B(v, e_i) = (1-\alpha)^k$ for any edge $e_i$ $(1 \leq i \leq m)$ in $E(\mathcal{S}_{m-3,0,1}^k)$. We also get $\sum_{e:v \in e}\left(B(v, e) + \alpha\right) = \rho_\alpha^\circ$, where $v$ is an arbitrary vertex of degree 1 in $\mathcal{S}_{m-3,0,1}^k$ or $v = u_2, u_3$.

Next, we compare $\sum_{e:u_1 \in e}\left(B(u_1, e) + \alpha\right)$ with $\rho_\alpha^\circ$. Since $C_2 > 0$ and $0 \leq \alpha < 1$, we have $\rho_\alpha^\circ - \alpha > \rho_\alpha^\circ - 2\alpha - C_2 > 0$ (by (42)). Therefore, we get $A_2 < \frac{(\rho_\alpha^\circ - \alpha)A_2}{\rho_\alpha^\circ - 2\alpha - C_2}$. It is noted that $\frac{(\rho_\alpha^\circ - \alpha)A_2}{\rho_\alpha^\circ - 2\alpha - C_2} > 0$. Furthermore, by (41), we have

$$\rho_\alpha^\circ - \sum_{e:u_1 \in e}\left(B(u_1, e) + \alpha\right)$$
$$= \rho_\alpha^\circ - \left(B(u_1, e_1) + (m-3)B(u_1, e_3) + (m-2)\alpha\right)$$
$$= -\alpha - A_2 + \frac{(\rho_\alpha^\circ - \alpha)A_2}{\rho_\alpha^\circ - 4\alpha - 3A_2} - \frac{(\rho_\alpha^\circ - \alpha)A_2}{\rho_\alpha^\circ - 2\alpha - C_2}$$

$$\leq -\alpha - 2A_2 + \frac{(\rho_\alpha^\circ - \alpha)A_2}{\rho_\alpha^\circ - 4\alpha - 3A_2}$$
$$= -\alpha + \frac{A_2}{\rho_\alpha^\circ - 4\alpha - 3A_2}\left(-\rho_\alpha^\circ + 7\alpha + 6A_2\right). \tag{43}$$

Bearing $\rho_\alpha^\circ > 4\alpha + 3A_2$ (by (39)), $\rho_\alpha^\circ > \alpha$, and $A_2 > 0$ in mind, when $m \geq 10$, by (41), we obtain

$$\rho_\alpha^\circ > (m-3)\alpha + (m-4)A_2 \geq 7\alpha + 6A_2. \tag{44}$$

Thus, when $m \geq 10$, by (43), we get $\rho_\alpha^\circ - \sum_{e:u_1 \in e}\left(B(u_1, e) + \alpha\right) < 0$. Therefore, $\mathcal{S}_{m-3,0,1}^k$ is strictly consistently $\rho_\alpha^\circ$-supernormal. By Lemma 3.3, we get $\rho_\alpha\left(\mathcal{S}_{m-3,0,1}^k\right) > \rho_\alpha^\circ = \rho_\alpha\left(\mathcal{S}_{3,m-4}^k\right)$.

(2). The proof of $\rho_\alpha\left(\mathcal{S}_{3,m-4}^k\right) > \rho_\alpha\bigl(\mathcal{T}(1, 2, m-4)\bigr)$, where $m = \frac{n-1}{k-1} \geq 10$.

In $\mathcal{T}(t_1, t_2, t_3)$ (as shown in Fig. 1(c)), let $t_1 = 1$, $t_2 = 2$ and $t_3 = m - 4$. Namely, we get $\mathcal{T}(1, 2, m-4)$. The vertices $u_1$, $u_2$, and $u_3$ and the edge $e_i$ ($1 \leq i \leq m$) of $\mathcal{T}(1, 2, m-4)$ are shown in Fig. 1(c). We construct a weighted incidence matrix $\boldsymbol{B}$ for $\mathcal{T}(1, 2, m-4)$ as follows. Let $B(v, e_i) = \rho_\alpha^\circ - \alpha$, where $v$ is an arbitrary vertex of degree 1 in $V\bigl(\mathcal{T}(1, 2, m-4)\bigr)$ and $e_i$ ($1 \leq i \leq m$) is the edge incident with $v$. Furthermore, let

$$B(u_1, e_1) = B_2, \qquad\qquad B(u_1, e_2) = A_2;$$
$$B(u_2, e_1) = \rho_\alpha^\circ - 3\alpha - 2A_2, \qquad B(u_2, e_3) = B(u_2, e_4) = A_2;$$
$$B(u_3, e_1) = \frac{(\rho_\alpha^\circ - \alpha)^2 A_2}{(\rho_\alpha^\circ - 3\alpha - 2A_2)B_2}, \qquad B(u_3, e_i) = A_2 \quad \text{for } 5 \leq i \leq m.$$

Since $A_2 > 0$ and $0 \leq \alpha < 1$, by (39), we have

$$\rho_\alpha^\circ - 3\alpha - 2A_2 > \rho_\alpha^\circ - 4\alpha - 3A_2 > 0. \tag{45}$$

Therefore, since $\rho_\alpha^\circ > \alpha$ and $A_2, B_2 > 0$, we get $B(v, e) > 0$ for any vertex $v$ and any edge $e$ incident with $v$ in $\mathcal{T}(1, 2, m-4)$.

It is easy to verify that $\prod_{v: v \in e_i} B(v, e_i) = (1-\alpha)^k$ for any edge $e_i$ ($1 \leq i \leq m$) in $E\bigl(\mathcal{T}(1, 2, m-4)\bigr)$. We can verify $\sum_{e:v \in e}\bigl(B(v,e) + \alpha\bigr) = \rho_\alpha^\circ$ for an arbitrary vertex $v$ of degree 1 in $\mathcal{T}(1, 2, m-4)$ or $v = u_1, u_2$.

Next, we compare $\sum_{e:u_3 \in e}\bigl(B(u_3, e) + \alpha\bigr)$ with $\rho_\alpha^\circ$. From (41), we have

$$\rho_\alpha^\circ - \sum_{e:u_3 \in e}\bigl(B(u_3, e) + \alpha\bigr)$$

$$
\begin{aligned}
&= \rho_\alpha^\circ - \big(B(u_3, e_1) + (m-4)B(u_3, e_5) + (m-3)\alpha\big) \\
&= \frac{(\rho_\alpha^\circ - \alpha)A_2}{\rho_\alpha^\circ - 4\alpha - 3A_2} - \frac{(\rho_\alpha^\circ - \alpha)^2 A_2}{(\rho_\alpha^\circ - 3\alpha - 2A_2)B_2} \\
&= D_2 \Big[(\rho_\alpha^\circ - 3\alpha - 2A_2) - \frac{(\rho_\alpha^\circ - \alpha)(\rho_\alpha^\circ - 4\alpha - 3A_2)}{B_2}\Big],
\end{aligned} \tag{46}
$$

where $D_2 = \frac{(\rho_\alpha^\circ - \alpha)A_2}{(\rho_\alpha^\circ - 4\alpha - 3A_2)(\rho_\alpha^\circ - 3\alpha - 2A_2)}$. It follows from $\rho_\alpha^\circ > \alpha$, $A_2 > 0$ and (45) that $D_2 > 0$. We can verify that $(\rho_\alpha^\circ - \alpha)(\rho_\alpha^\circ - 4\alpha - 3A_2) - B_2(\rho_\alpha^\circ - 3\alpha - 2A_2) = -2(\alpha + A_2)^2 < 0$. Therefore, it follows from (46) that $\rho_\alpha^\circ > \sum_{e:u_3 \in e}\big(B(u_3, e) + \alpha\big)$. Therefore, $\mathcal{T}(1, 2, m-4)$ is strictly $\rho_\alpha^\circ$-subnormal. By Lemma 3.2, we have $\rho_\alpha\big(\mathcal{S}_{3,m-4}^k\big) > \rho_\alpha\big(\mathcal{T}(1, 2, m-4)\big)$. □

**Lemma 4.3** Let $k \geq 3$ and $m = \frac{n-1}{k-1} \geq 13$. We have $\rho_\alpha\big(\mathcal{T}(1, 2, m-4)\big) > \rho_\alpha\big(\mathcal{S}_{4,m-5}^k\big)$.

**Proof**. In $\mathcal{S}_{a,b}^k$ (as shown in Fig. 1(a)), let $a = 4$ and $b = m-5$. Namely, we get $\mathcal{S}_{4,m-5}^k$. For simplicity, let $\rho_\alpha^\Delta = \rho_\alpha\big(\mathcal{S}_{4,m-5}^k\big)$. Let $\boldsymbol{x} = (x_1, \ldots, x_n) \in \mathbb{R}_{++}^n$ be the $\alpha$-Perron vector of $\rho_\alpha^\Delta$. We suppose that $\mathcal{S}_{4,m-5}^k$ is consistently $\rho_\alpha^\Delta$-normal. Let $v_1$, $v_2$ and $v_3$ of $\mathcal{S}_{4,m-5}^k$ be three vertices of degree 1 which are respectively incident with $e_1$, $e_2$ ad $e_m$ of $\mathcal{S}_{4,m-5}^k$, where the vertices $v_1$, $v_2$, $v_3$, $u_1$, and $u_2$ and the edges $e_i$ ($1 \leq i \leq m$) of $\mathcal{S}_{4,m-5}^k$ are shown in Fig. 1(a). By the eigenequations (4) of $\mathcal{S}_{4,m-5}^k$ at $v_1$, $v_2$, $v_3$, $u_1$, and $u_2$ and bearing the symmetry of the entries in $\boldsymbol{x}$ in mind, we get

$$\rho_\alpha^\Delta x_{v_1}^{k-1} = \alpha x_{v_1}^{k-1} + (1-\alpha)x_{u_1}x_{u_2}x_{v_1}^{k-3}, \tag{47}$$

$$\rho_\alpha^\Delta x_{v_2}^{k-1} = \alpha x_{v_2}^{k-1} + (1-\alpha)x_{u_1}x_{v_2}^{k-2}, \tag{48}$$

$$\rho_\alpha^\Delta x_{v_3}^{k-1} = \alpha x_{v_3}^{k-1} + (1-\alpha)x_{u_2}x_{v_3}^{k-2}, \tag{49}$$

$$\rho_\alpha^\Delta x_{u_1}^{k-1} = 5\alpha x_{u_1}^{k-1} + 4(1-\alpha)x_{v_2}^{k-1} + (1-\alpha)x_{u_2}x_{v_1}^{k-2}, \tag{50}$$

$$\rho_\alpha^\Delta x_{u_2}^{k-1} = (m-4)\alpha x_{u_2}^{k-1} + (m-5)(1-\alpha)x_{v_3}^{k-1} + (1-\alpha)x_{u_1}x_{v_1}^{k-2}. \tag{51}$$

By (47)–(49), we get respectively

$$x_{v_1} = \sqrt{\frac{1-\alpha}{\rho_\alpha^\Delta - \alpha} x_{u_1} x_{u_2}}, \quad x_{v_2} = \frac{1-\alpha}{\rho_\alpha^\Delta - \alpha} x_{u_1}, \quad x_{v_3} = \frac{1-\alpha}{\rho_\alpha^\Delta - \alpha} x_{u_2}. \tag{52}$$

From (47), we have $\rho_\alpha^\Delta > \alpha$ when $\boldsymbol{x} \in \mathbb{R}_{++}^n$ and $0 \leq \alpha < 1$. For simplicity, let

$$A_3 = \frac{(1-\alpha)^k}{(\rho_\alpha^\Delta - \alpha)^{k-1}}, \quad B_3 = \rho_\alpha^\Delta - 2\alpha - A_3. \tag{53}$$

Since $\rho_\alpha^\Delta > \alpha$ and $0 \leq \alpha < 1$, we have $A_3 > 0$.

Substituting (52) into (50), we get

$$\rho_\alpha^\Delta > 5\alpha + 4A_3, \quad (\text{since } \boldsymbol{x} \in \mathbb{R}_{++}^n \text{ and } 0 \leq \alpha < 1), \tag{54}$$

$$x_{u_1} = \frac{1-\alpha}{(\rho_\alpha^\Delta - \alpha)^{1-2/k}(\rho_\alpha^\Delta - 5\alpha - 4A_3)^{2/k}} x_{u_2}. \tag{55}$$

Substituting (52) and (55) into (51), we obtain

$$\rho_\alpha^\Delta - (m-4)\alpha = (m-5)A_3 + \frac{(\rho_\alpha^\Delta - \alpha)A_3}{\rho_\alpha^\Delta - 5\alpha - 4A_3}. \tag{56}$$

We construct a weighted incidence matrix $\boldsymbol{B}$ for $\mathcal{T}(1, 2, m-4)$ as follows. The vertices $u_1$, $u_2$, and $u_3$ and the edge $e_i$ ($1 \leq i \leq m$) of $\mathcal{T}(1, 2, m-4)$ are shown in Fig. 1(c) with $t_1 = 1$, $t_2 = 2$ and $t_3 = m-4$. Let $B(v, e_i) = \rho_\alpha^\Delta - \alpha$, where $v$ is an arbitrary vertex of degree 1 in $V(\mathcal{T}(1, 2, m-4))$ and $e_i$ is the edge incident with $v$. Furthermore, let

$$B(u_1, e_1) = B_3, \quad B(u_2, e_1) = \rho_\alpha^\Delta - 3\alpha - 2A_3;$$

$$B(u_1, e_2) = B(u_2, e_3) = B(u_2, e_4) = B(u_3, e_i) = A_3 \quad \text{for } 5 \leq i \leq m;$$

$$B(u_3, e_1) = \frac{(\rho_\alpha^\Delta - \alpha)^2 A_3}{(\rho_\alpha^\Delta - 3\alpha - 2A_3)B_3}.$$

Since $A_3 > 0$ and $0 \leq \alpha < 1$, by (53) and (54), we obtain $B_3 > 0$ and

$$\rho_\alpha^\Delta - 3\alpha - 2A_3 > \rho_\alpha^\Delta - 5\alpha - 4A_3 > 0. \tag{57}$$

Furthermore, since $\rho_\alpha^\Delta > \alpha$ and $A_3, B_3 > 0$, we obtain $B(v, e) > 0$ for any vertex $v$ and any edge $e$ incident with $v$ in $\mathcal{T}(1, 2, m-4)$.

It is easy to verify that $\prod_{v:v\in e_i} B(v, e_i) = (1-\alpha)^k$ for any edge $e_i$ ($1 \leq i \leq m$) in $E(\mathcal{T}(1, 2, m-4))$. We can also check that $\sum_{e:v\in e}(B(v,e) + \alpha) = \rho_\alpha^\Delta$ for an arbitrary vertex $v$ of degree 1 in $V(\mathcal{T}(1, 2, m-4))$ or $v = u_1, u_2$.

Next, we compare $\sum_{e:u_3 \in e}(B(u_3, e) + \alpha)$ with $\rho_\alpha^\Delta$. Since $\rho_\alpha^\Delta - 3\alpha - 2A_3 > 0$ (by (57)) and $A_3, B_3 > 0$, we have $\frac{(\rho_\alpha^\Delta - \alpha)^2 A_3}{(\rho_\alpha^\Delta - 3\alpha - 2A_3)B_3} > 0$. Since $A_3 > 0$ and $0 \leq \alpha < 1$, we have $(\rho_\alpha^\Delta - \alpha)^2 > (\rho_\alpha^\Delta - 3\alpha - 2A_3)(\rho_\alpha^\Delta - 2\alpha - A_3) = (\rho_\alpha^\Delta - 3\alpha - 2A_3)B_3$. Therefore, we obtain $A_3 < \frac{(\rho_\alpha^\Delta - \alpha)^2 A_3}{(\rho_\alpha^\Delta - 3\alpha - 2A_3)B_3}$. Furthermore, by (56), we have

$$\rho_\alpha^\Delta - \sum_{e:u_3 \in e}(B(u_3, e) + \alpha)$$

$$= \rho_\alpha^\Delta - \big(B(u_3, e_1) + (m-4)B(u_3, e_5) + (m-3)\alpha\big)$$

$$= -\alpha - A_3 + \frac{(\rho_\alpha^\Delta - \alpha)A_3}{\rho_\alpha^\Delta - 5\alpha - 4A_3} - \frac{(\rho_\alpha^\Delta - \alpha)^2 A_3}{(\rho_\alpha^\Delta - 3\alpha - 2A_3)B_3}$$

$$\leq -\alpha - 2A_3 + \frac{(\rho_\alpha^\Delta - \alpha)A_3}{\rho_\alpha^\Delta - 5\alpha - 4A_3}$$

$$= -\alpha + \frac{A_3}{\rho_\alpha^\Delta - 5\alpha - 4A_3}\left(-\rho_\alpha^\Delta + 9\alpha + 8A_3\right). \tag{58}$$

Since $\rho_\alpha^\Delta > 5\alpha + 4A_3$ (by (54)), $\rho_\alpha^\Delta > \alpha$, and $A_3 > 0$, it follows from (56) that $\rho_\alpha^\Delta > (m-4)\alpha + (m-5)A_3 \geq 9\alpha + 8A_3$ when $m \geq 13$. Thus, by (58), we get $\rho_\alpha^\Delta - \sum_{e:u_3 \in e}\left(B(u_3, e) + \alpha\right) < 0$ when $m \geq 13$. Therefore $\mathcal{S}_{m-3,0,1}^k$ is strictly consistently $\rho_\alpha$-supernormal. By Lemma 3.3, we get $\rho_\alpha\bigl(\mathcal{T}(1, 2, m-4)\bigr) > \rho_\alpha^\Delta = \rho_\alpha\bigl(\mathcal{S}_{4,m-5}^k\bigr)$. □

By Lemmas 2.6–2.8 and 4.1–4.3, we finally obtain that Conjecture 1 is correct for $0 \leq \alpha < 1$, which is shown in Theorem 4.1.

**Theorem 4.1** *Let $0 \leq \alpha < 1$, $k \geq 3$ and $m = \frac{n-1}{k-1} \geq 13$. For any supertree $\mathcal{T} \in \mathcal{T}(n, k) \backslash \mathcal{T}^{(8)}(n, k)$, we have*

$$\rho_\alpha\left(\mathcal{S}_{m+1}^k\right) > \rho_\alpha\left(\mathcal{S}_{1,m-2}^k\right) > \rho_\alpha\left(\mathcal{S}_{2,m-3}^k\right) > \rho_\alpha\left(\mathcal{T}(1,1,m-3)\right) >$$
$$\rho_\alpha\left(\mathcal{S}_{1,m-4,1}^k\right) > \rho_\alpha\left(\mathcal{S}_{m-3,0,1}^k\right) > \rho_\alpha\left(\mathcal{S}_{3,m-4}^k\right) > \rho_\alpha\left(\mathcal{T}(1,2,m-4)\right) > \rho_\alpha\left(\mathcal{T}\right).$$

**Proof**. Let $k \geq 3$ and $m = \frac{n-1}{k-1} \geq 13$. You et al. [12] obtained $\rho_\alpha\left(\mathcal{S}_{m+1}^k\right) > \rho_\alpha\left(\mathcal{S}_{1,m-2}^k\right) > \rho_\alpha\left(\mathcal{S}_{2,m-3}^k\right) > \rho_\alpha\left(\mathcal{T}(1,1,m-3)\right)$. By Lemmas 4.1 and 4.2, we get $\rho_\alpha\left(\mathcal{T}(1,1,m-3)\right) > \rho_\alpha\left(\mathcal{S}_{1,m-4,1}^k\right) > \rho_\alpha\left(\mathcal{S}_{m-3,0,1}^k\right) > \rho_\alpha\left(\mathcal{S}_{3,m-4}^k\right) > \rho_\alpha\left(\mathcal{T}(1,2,m-4)\right)$.

Let $\mathcal{T}$ be an arbitrary supertree in $\mathcal{T}(n, k) \backslash \mathcal{T}^{(8)}(n, k)$. If $\mathcal{T}$ has only two vertices which have degrees greater than one, then $\mathcal{T} \in \{\mathcal{S}_{4,m-5}^r, \mathcal{S}_{5,m-6}^r, \ldots \mathcal{S}_{\lfloor\frac{m-1}{2}\rfloor, \lceil\frac{m-1}{2}\rceil}^r\}$. By Lemmas 4.3 and 2.8, we get $\rho_\alpha\bigl(\mathcal{T}(1,2,m-4)\bigr) > \rho_\alpha\left(\mathcal{S}_{4,m-5}^k\right) > \rho_\alpha\left(\mathcal{T}\right)$. Next, we assume that $\mathcal{T}$ has at least three vertices which have degrees greater than one. Let $\pi = (d_0, d_1, \ldots, d_{n-1})$ be the degree sequence of $\mathcal{T}$, where $2 \leq d_2 \leq d_1 \leq d_0 \leq m-3$. Obviously, the degree sequence $\pi$ of $\mathcal{T}$ is different from the degree sequence of $\mathcal{S}_{j,m-j}^r$, where $0 \leq j \leq \lfloor\frac{m-1}{2}\rfloor$. Let $\pi_3' = (m-3, 3, 2, 1, \ldots, 1)$. Obviously, the degree sequence of $\mathcal{T}(1, 2, m-4)$ is $\pi_3'$. Bearing the definition of $\mathcal{T}(1, 2, m-4)$ in mind, we get that $\mathcal{T}(1, 2, m-4)$ has a BFS-ordering. Thus, by Lemma 2.6, $\mathcal{T}(1, 2, m-4)$ is the unique $k$-uniform supertree with the largest $\alpha$-spectral radius in $T_{\pi_3'}$. Therefore, if $\pi = \pi_3'$, then $\rho_\alpha\bigl(\mathcal{T}(1,2,m-4)\bigr) > \rho_\alpha\left(\mathcal{T}\right)$. If $\pi \neq \pi_3'$, we have $\pi \lhd \pi_3'$. Thus, by Lemma 2.7, we get $\rho_\alpha\bigl(\mathcal{T}(1,2,m-4)\bigr) > \rho_\alpha\left(\mathcal{T}\right)$. By combining the above proofs, we obtain Theorem 4.1. □

## Acknowledgments

The work was supported by the Natural Science Foundation of Shanghai under the grant number 21ZR1423500.

---

[1] J. Cooper, A. Dutle, Spectra of uniform hypergraphs, Linear Algebra and its Applications 436 (2012) 3268–3292.

[2] V. Nikiforov, Merging the A- and Q-spectral theories, Applicable Analysis and Discrete Mathematics 11 (2017) 81–107.

[3] H. Y. Lin, H. Y. Guo, B. Zhou, On the $\alpha$-spectral radius of irregular uniform hypergraphs, Linear and Multilinear Algebra 68 (2020) 265–277.

[4] F. F. Wang, H. Y. Shan, Z. Y. Wang, On some properties of the $\alpha$-spectral radius of the $k$-uniform hypergraph, Linear Algebra and its Applications 589 (2020) 62–79.

[5] H. Y. Lin, B. Zhou, The $\alpha$-spectral radius of general hypergraphs, Applied Mathematics and Computation 386 (2020) 125449.

[6] L. Y. Kang, J. Wang, E. F. Shan, The maximum $\alpha$-spectral radius of unicyclic hypergraphs with fixed diameter, Acta Mathematica Sinica, English Series 38 (2022) 924–936.

[7] Y. Hou, A. Chang, C. Shi, On the $\alpha$-spectra of uniform hypergraphs and its associated graphs, Acta Mathematica Sinica, English Series 36 (2020) 842–850.

[8] H. Y. Guo, B. Zhou, On the $\alpha$-spectral radius of uniform hypergraphs, Discussiones Mathematicae Graph Theory 40 (2020) 559–575.

[9] C. X. Duan, L. G. Wang, The $\alpha$-spectral radius of $f$-connected general hypergraphs, Applied Mathematics and Computation 382 (2020) 125336.

[10] H. H. Li, J. Y. Shao, L. Q. Qi, The extremal spectral radii of $k$-uniform supertrees, Journal of Combinatorial Optimization 32 (2016) 741–764.

[11] X. Y. Yuan, J. Y. Shao, H. Y. Shan, Ordering of some uniform supertrees with larger spectral radii, Linear Algebra and its Applications 495 (2016) 206–222.

[12] L. H. You, L. H. Deng, Y. F. Huang, The maximum $\alpha$-spectral radius and the majorization theorem of $k$-uniform supertrees, Discrete Applied Mathematics 285 (2020) 663–675.

[13] W. H. Wang, J. X. Zhou, R. Sun, On the conjecture of the r-uniform supertrees with the

eight largest -spectral radii, Discrete Applied Mathematics 322 (2022) 311–319.

[14] L. Y. Lu, S. D. Man, Connected hypergraphs with small spectral radius, Linear Algebra and its Applications 509 (2016) 206–227.

[15] S. Friedland, S. Gaubert, L. Han, Perron–Frobenius theorem for nonnegative multilinear forms and extensions, Linear Algebra and its Applications 438 (2013) 738–749.

[16] Y. N. Yang, Q. Z. Yang, On some properties of nonnegative weakly irreducible tensors, arXiv:1111.0713v3 (2011).

[17] Y. N. Yang, Q. Z. Yang, Further results for Perron–Frobenius theorem for nonnegative tensors, SIAM Journal on Matrix Analysis and Applications 31 (2010) 2517–2530.

[18] K. J. Pearson, T. Zhang, On spectral hypergraph theory of the adjacency tensor, Graphs and Combinatorics 30 (2014) 1233–1248.

[19] L. Q. Qi, Symmetric nonnegative tensors and copositive tensors, Linear Algebra and its Applications 439 (2013) 228–238.

[20] P. Xiao, L. G. Wang, Y. Lu, The maximum spectral radii of uniform supertrees with given degree sequences, Linear Algebra and its Applications 523 (2017) 33–45.